\documentclass[12pt]{article}
\usepackage[all]{xy}
\usepackage{amsmath}
\usepackage{amsfonts}
\usepackage{amssymb}
\usepackage{amscd}
\usepackage{amsthm}
\usepackage{latexsym}
\usepackage{amsbsy}
\theoremstyle{plain}
\newtheorem{theorem}{Theorem}[section]
\newtheorem{lemma}[theorem]{Lemma}
\newtheorem{proposition}[theorem]{Proposition}
\newtheorem{corollary}[theorem]{Corollary}
\newtheorem{definition}[theorem]{Definition}
\newtheorem{example}[theorem]{Example}

\newcommand{\ot}{\otimes}
\date{}
\begin{document}
\title
{Invariant Cyclic Homology}

\author {M. Khalkhali,~~~ B. Rangipour,
\\\texttt{~masoud@uwo.ca ~~~~brangipo@uwo.ca}
\\ Department of Mathematics 
\\ University of Western Ontario 
\\ London ON, Canada}
\maketitle
\begin{abstract}
We define a noncommutative analogue of invariant de Rham cohomology.
 More precisely,  for a triple $(A,\mathcal{H},M)$ consisting of a Hopf algebra $\mathcal{H}$, an
  $\mathcal{H}$-comodule 
 algebra $A$,   an $\mathcal{H}$-module $M$, and  a compatible grouplike element $\sigma$ in $\mathcal{H}$, we define 
  the cyclic module of invariant chains on $A$ with coefficients in $M$ and call its cyclic homology the invariant
   cyclic homology of $A$ with coefficients in $M$. We also develop a dual theory for
coalgebras. Examples include cyclic cohomology of Hopf algebras defined by Connes-Moscovici and its dual theory.
 We establish various results and computations including  one for the quantum group $SL(q,2)$.
\end{abstract}
\section{Introduction}
 It is well known that cyclic homology replaces de Rham cohomology in noncommutative settings.
 For example, by a result of Connes~\cite{ac85} the periodic cyclic homology of the algebra of smooth
 functions on a smooth manifold  is isomorphic to the de Rham cohomology of the manifold. Invariant de Rham 
 cohomology  was introduced by Chevalley  and Eilenberg in~\cite{ce} in part to
 relate the de Rham  cohomology of a Lie group to the cohomology of its Lie algebra. 
   
 In this paper we define a  noncommutative analogue of invariant  de Rham cohomology. 
 More precisely, for a triple $(A,\mathcal{H},M)$ consisting of a Hopf algebra $\mathcal{H}$,
 an $\mathcal{H}$\nobreakdash-comodule algebra $A$ and an $\mathcal{H}$\nobreakdash-module $M$, and for a suitably 
 chosen grouplike element $\sigma\in \mathcal{H}$,  we define the 
  {\it cyclic module of invariant chains}  on $A$, denoted $\{C_n^\mathcal{H}(A,M)\}_n$,  as the space of
   coinvariants
   of the paracyclic
 module $\{C_n(A,M)\}_n$ under the action of $\mathcal{H}$. Such triples are called \mbox{\it $\sigma$-compatible
  Hopf
  triple} in this paper. We call the resulting cyclic homology groups, the {\it invariant cyclic homology 
 of the algebra} $A$ with coefficients  in $M$. For example  in complete analogy with the classical   case,
   the  cyclic module $\{\widetilde{\mathcal{H}}_n^{(\delta,\sigma)}\}_{n}$, defined in ~\cite{kr2} for a 
     Hopf algebra $\mathcal{H}$ endowed  with a modular pair $(\delta, \sigma)$ in involution, is isomorphic with
      the 
     invariant cyclic module 
  of the algebra $\mathcal{H}$   
  with respect to the coaction of $\mathcal{H}$ on itself via the
  comultiplication map $\mathcal{H}\rightarrow \mathcal{H}\ot\mathcal{H}$. (Classically ~\cite{ce}, the 
  cohomology of the Lie 
  algebra of a Lie group $G$ appeared  as the invariant de Rham cohomology of $G$ with respect to the 
  action of $G$ on itself via multiplication).
        
 We also develop a  theory for coalgebras. This is needed, for instance, 
  in order to treat the Connes-Moscovici cocyclic
 module $\{\mathcal{H}^n_{(\delta, \sigma)}\}_{n}$  of a Hopf algebra $\mathcal{H}$ as the 
 invariant cocyclic module of $\mathcal{H}$; but we go beyond this case and introduce
 $\delta$\nobreakdash-{\it compatible Hopf cotriples } $(C,\mathcal{H},V)$.
  Thus for an $\mathcal{H}$\nobreakdash-module coalgebra $C$, an  
  $\mathcal{H}$\nobreakdash-comodule $V$, and a compatible character $\delta$ 
  we define the {\it  cocyclic module of invariant cochains} on $C$, 
   $\{C_\mathcal{H}^n(C,V)\}_n$,  as the space
 of coinvariants of the paracocyclic module  $\{C^n(C, V)\}_n$. 
 
 One feature of our approach is that even for $A=\mathcal{H}$ or $C=\mathcal{H}$ we allow coefficients  to 
 enter the picture and thus can interpret  group homology and Lie algebra homology 
  with coefficients as invariant cyclic homology.  

There  are two other situations, corresponding to Hopf module algebras and Hopf 
comoule coalgebras, where one can also define an invariant cyclic homology theory. 
In particular, the twisted cyclic cohomology of \cite{kmt} is a special case of 
invariant cyclic cohomology theory corresponding to $G=\mathbb{Z}$. There are also
 interesting connections with covariant differential calculi on Hopf algebras and 
quantum groups. These  matters will be pursued elsewhere \cite{kr6}.

   This work grew out of our attempt to interpret the Connes-Moscovici  cocyclic module for Hopf algebras as well as the cyclic module 
   in {\cite{kr2,tr}} as special case of a general invariant cyclic homology theory.   
   The existence of invariant cyclic homology and the methods we use was inspired by their work. We  also 
   mention that  M. Crainic in  \cite{cr} has interpreted    the Connes-Moscovici cocyclic module $\{\mathcal{H}_{(\delta,1)}^n\}_n$ as 
      the space of invariant cochains on $\mathcal{H}$.       

 \section{Preliminaries}
   Let $k$ be a commutative unital ring. In this   paper, by an algebra we mean a unital  associative 
  algebra over $k$.
    The same convention applies to coalgebras and Hopf algebras. We denote the coproduct of a coalgebra 
  by $\Delta$,
 its counit by $\epsilon$, and the antipode of a Hopf algebra  by $S$. The unadorned tensor 
  product $\ot$ means tensor product over $k$. We use Sweedler's notation and through this paper write 
  $\Delta(h)=h^{(1)}\otimes h^{(2)}$,
  $\Delta^2(h):=(\Delta\ot id)\circ\Delta(h)=h^{(1)}\ot h^{(2)}\ot h^{(3)}$, where  summation is understood.
  A {\it character}  of a Hopf algebra $\mathcal{H}$ is a (unital ) algebra homomorphism 
  $\delta:\mathcal{H}\longrightarrow k$. A {\it grouplike} element is an element $\sigma\in \mathcal{H}$
  such that $\Delta(\sigma)=\sigma\ot \sigma$  and $\epsilon(\sigma)=1$.
             
 Let $\mathcal{H}$ be a Hopf algebra. By an $\mathcal{H}$\nobreakdash-module we
  mean a module over the underlying algebra  of $\mathcal{H}$ and by an $\mathcal{H}$\nobreakdash-comodule 
 we mean a comodule over the underlying coalgebra of $\mathcal{H}$. We use Sweedler's notation for comodules.
 Thus if $\rho:  M\longrightarrow C\ot M$ is the structure map of a left $C$\nobreakdash-comodule $M$, 
 we write $\rho(m)=m^{(-1)}\ot m^{(0)}$, where  summation is understood.
               
 Let $\mathcal{H}$ be a Hopf algebra. An algebra $A$ is called a left $\mathcal{H}$-{\it comodule algebra}
 if $A$ is a  left  $\mathcal{H}$\nobreakdash-comodule via   $\rho:  A\longrightarrow  \mathcal{H}\ot A$ and $\rho$  
 is an algebra map.                
 Similarly, a coalgebra $C$ is called a left  $\mathcal{H}$-{\it module coalgebra}
 if  $C$ is a  left  $\mathcal{H}$\nobreakdash-module via   $\mu:  \mathcal{H}\ot C\longrightarrow C$ and $\mu$  
 is a coalgebra map.
                
 Let $M$ be a left $\mathcal{H}$\nobreakdash-module and $\delta $ a character of $\mathcal{H}$.
 The space of {\it coinvariants}  of $M$ (with respect to $\delta$) is the  $k$\nobreakdash-module
 $$M_\mathcal{H}=M/ \text{span}\{ hm-\delta(h)m\mid  h\in\mathcal{H},\; m\in             
      M\}=k_\delta\underset{\mathcal{H}}{\ot}M,$$
  where the right $\mathcal{H}$\nobreakdash-module $k_\delta$ is defined by $k_\delta=k$ and $\mathcal{H}$ acts via the character 
  $\delta$. 
 Similarly, if  $M$ is a left  $\mathcal{H}$\nobreakdash-comodule, and $\sigma\in \mathcal{H}$ 
 is a grouplike element, the space of  coinvariants  of $M$ is the $k$\nobreakdash-module 
 $$M^{\text{co}\mathcal{H}}=\{ m\in M\mid \rho(m)= \sigma\ot m \}=Hom_\mathcal{H}(k_\sigma,M),$$ 
 where $k_\sigma =k$ is the left $\mathcal{H}$\nobreakdash-comodule defined by   $\sigma$.
                 
  These concepts  are usually considered for $\delta=\epsilon$, the counit of $\mathcal{H}$ and 
  $\sigma=1$, the unit of $\mathcal{H}$. The work of Connes-Moscovici ~\cite{achm98,achm99,achm00}, however, 
  shows that for cyclic cohomology 
  of Non(co)commutative Hopf algebras, it is absolutely necessary to consider this more general case.   
                  
 By a {\it paracyclic} module we mean a simplicial module $M=\{M_n\}_{n}$ endowed with $k$-linear  maps 
 $\tau_n: M_n\longrightarrow M_n$, $n\ge 0$, such that the following identities hold                  
 \begin{eqnarray*}
{\delta }_i {\tau}_n &=& {\tau}_{n-1}{\delta }_{i-1} \\
\delta _0 {\tau}_n & =& {\delta }_n \\
{\sigma}_i {\tau}_n &=& {\tau}_{n+1}{\sigma}_{i-1} \\
{\sigma}_0 {\tau}_n &=& {\tau}_{n+1}^2{\sigma}_{n}. 
\end{eqnarray*}    
  
Here $\delta_i$ and $\sigma_i$ denote the face and the degeneracy operators of $M$. If $\tau^{n+1}_n=id$, for all 
$n\ge 0$, we have a {\it cyclic module} in the sense of Connes~\cite{ac85}. For example, 
if $A$ is an algebra and $g:A \longrightarrow A$   is an automorphism of $A$, 
one can check that the following operators define a paracyclic module denoted by  $A_g^\natural$,
 where $A^\natural_{g,n}=A^{\ot(n+1)}$ and 
 \begin{eqnarray*}
 &&\delta_0(a_0\ot a_1\ot\dots\ot a_n)=a_0a_1\ot a_2\ot\dots \ot a_n,\\
 &&\delta_i(a_0\ot a_1\ot\dots\ot a_n)=a_0\ot \dots \ot a_i a_{i+1}\ot\dots \ot  a_n, \qquad 1\le i\le n-1,\\
  &&\delta_n(a_0\ot a_1\ot\dots\ot a_n)=g(a_n)a_0\ot a_1\dots \ot a_{n-1},\\        
&&\sigma_i(a_0\ot a_1\ot\dots\ot a_n)=a_0 \ot \dots \ot a_i\ot 1 \ot \dots \ot a_n, \qquad 0\le i\le n-1,\\  
&&\tau(a_0\ot a_1\ot\dots\ot a_n)=g(a_n)\ot a_0\ot a_1\dots \ot a_{n-1}. 
\end{eqnarray*}
Since $\tau^{n+1}(a_0\ot\dots\ot a_n)=ga_0\ot\dots \ot ga_n$, it is clear that $A_g^\natural$ is a cyclic module 
if and only if $g=id$.

Dually, we have the notion  of a {\it paracocyclic module}. For example if $C$ is a coalgebra 
and $\theta: C\longrightarrow C $  an automorphism of $C$, the paracocyclic module $C_\natural^\theta$
 is defined by $C_\natural^{\theta,n}=C^{\ot(n+1)},\; n\ge 0$,  and  
\begin{eqnarray*}
&&\delta_i(c_0 \otimes c_1 \otimes \dots \otimes c_n)= c_0 \otimes\dots 
\otimes  c_i^{(1)}\otimes c_i^{(2)}\otimes c_n ~~~ 0 \leq i \leq n\\
&&\delta_{n+1}(c_0 \otimes c_1 \otimes \dots \otimes c_n)= c_0^{(2)}\otimes c_1 
\otimes \dots \otimes c_n  \otimes \theta(c_0^{(1)})\\
&&\sigma_i(c_0 \otimes c_1 \otimes \dots \otimes c_n)=c_0 \otimes \dots 
c_i \otimes \varepsilon(c_{i+1})\otimes 
\dots\otimes c_n ~~~0\leq i \leq n-1\\
&&\tau(c_0 \otimes c_1 \otimes \dots \otimes c_n)  = c_1 \otimes c_2 \otimes 
\dots \otimes c_n \otimes \theta(c_0).
\end{eqnarray*}

  We denote the cyclic (co)homology group of a (co)cyclic module
   $M$ by $HC_\bullet(M)$ (resp. $HC^\bullet(M))$ and if 
   $M=A^\natural$ or $M=C_\natural$, we denote them by $HC_\bullet(A)$ and $HC^\bullet(C)$ respectively. 
      
          We recall the statement of the Eilenberg-Zilber theorem for cylindrical 
          modules from \cite{gj}(cf. \cite{kr1} for a purely algebraic proof). This result
                     is needed in Subsection \ref{subsection3.4}. Recall  that a {\it cylindrical module } 
                     is a doubly  graded $k$\nobreakdash-module
            $\{ X_{p,q}\}_{p,q}$ such that each row and each column is a 
            paracyclic module, all horizontal operators commute with all vertical  
            operators and \mbox{for all $p,q$,  ${\tau^{p+1}t^{q+1}=id :X_{p,q}\longrightarrow X_{p,q}}$.} Here  
           \mbox{ ${t,\tau: X_{p,q}\longrightarrow X_{p,q}}$} denote the horizontal and vertical cyclic operators     
            respectively. The horizontal simplicial operators are denoted by $d_i$, $s_i$, and 
            the vertical operators by $\delta_i$, $\sigma_i$.
            
            Given a cylindrical module $X$, its {\it diagonal}, denoted by $d(X)$ is defined by 
            $d(X)=X_{n,n}$ with simplicial and cyclic operators given by $d_i\delta_i$,
             $s_i\sigma_i$ and $t\tau$. It is a cyclic module. 
            Associated to this cyclic module, we have a mixed complex denoted by $(d(X),b_d,B_d)$. 
            The {\it Total complex}
            of $X$ is also a mixed complex  $(Tot(X),b_t,B_t)$ with 
            \mbox{$Tot(X)_n=\underset{p+q=n}{\bigoplus}X_{p,q}$.}
             The { generalized cyclic Eilenberg-Zilber} theorem \cite{gj} 
             states that these two mixed complexes are 
             chain homotopy equivalent. In particular they have isomorphic cyclic homology groups.  
             
             Finally we make an important remark about the antipode of Hopf algebras in this paper.
               Throughout this paper we assume that the antipode is bijective. We 
             need this hypothesis throughout this paper and in particular in the proof of
              Theorems \ref{1} and \ref{theorem4.10}. It is known that the class  
               of Hopf algebras with bijective antipode is a very large class that contains quantum groups.

 \section{Invariant Cyclic Homology of Hopf Triples}\label{sec3}    
 
 In this section we define the concept of $\sigma$\nobreakdash-{\it compatible Hopf triple} and its invariant cyclic
  homology. In Subsection \ref{subsection3.1} we give  various examples and calculations including one 
  involving the quantum group $A(SL_q(2))$ and a $2$-dimensional module.  We also prove a Morita invariance theorem.

  In Subsection \ref{subsection3.2}  we compare the invariant cyclic homology with the invariant
   de Rham cohomology for smooth affine algebras  . 
  The comparison maps should be isomorphisms  but we can verify this only for algebraic groups.

  In Subsection {\ref{subsection3.4}} we consider the invariant cyclic homology of smash 
  product algebras $A=\mathcal{H}\#B$, where $B$ is an $\mathcal{H}$-module algebra. We derive a spectral sequence 
  and show that it collapses if $\mathcal{H}$ is semisimple. The ultimate goal here would be
  to establish a similar spectral  sequence for all Hopf-Galois extensions, but
   it is not clear how to do this at the moment.

 \subsection{Definition and basic properties}\label{subsection3.1}
 \begin{definition}
 By a left  {\it Hopf triple } we mean a  triple $(A, \mathcal{H}, M)$, where $\mathcal{H}$ is a Hopf algebra,
  $A$ is a left $\mathcal{H}$\nobreakdash-comodule algebra and $M$ is a left $\mathcal{H}$\nobreakdash-module.
  Right Hopf triples are defined in a similar way.     
\end{definition}  
\begin{example}\label{example3.2}~~~
\begin{itemize}
\item [(i)](Trivial triples). Let $\mathcal{H}=k$, $M$ any $k$\nobreakdash-module, and $A$ any $k$-algebra. 
Then $(A,k,M)$ is a left Hopf triple.
\item[(ii)] Let $\mathcal{H}$ be a Hopf algebra and $M$ a left $\mathcal{H}$\nobreakdash-module.  Then 
$(A,\mathcal{H},M)$ is a left Hopf triple,  where $A=\mathcal{H}$ is the underlying algebra of $\mathcal{H}$
 and $\mathcal{H}$ coacts on $\mathcal{H}$ via its comultiplication. In particular, for $M=k$ and $\mathcal{H}$     
 acting on $k$ via a character $\delta$, we obtain a Hopf triple $(\mathcal{H},\mathcal{H},k_\delta)$.
  \item[(iii)]  Let $G$ be an affine algebraic group acting from left 
  on an affine algebraic variety $X$.  Let $\mathcal{H}=\mathbb{C}\lbrack G\rbrack$ and 
  $A=\mathbb{C}\lbrack X\rbrack$
   be the coordinate rings of $G$ and $X$, respectively.
    Then $A$ is a left  $\mathcal{H}$\nobreakdash-comodule algebra.
    The coaction
    \mbox{$\rho: A\rightarrow \mathcal{H}\otimes A$} is induced by the action \mbox{$G\times X\longrightarrow X$}. 
    Thus we obtain a Hopf triple \mbox{$(\mathbb{C}\lbrack X\rbrack,\mathbb{C}\lbrack G\rbrack, \mathbb{C}).$ } 
 \item[(iv)] Let $G$ be a group. It is easy to see that a $kG$\nobreakdash-comodule algebra is nothing but a $G$-graded
  algebra $A=\underset{g\in G}{\oplus}A_g$ ~\cite{mo}. The coaction  $\rho : A\longrightarrow kG\ot A$
  is given by $\rho(a)=\underset{g}{\sum}g\ot a_g$, where $a=\underset{g}{\sum}a_g$.
   Thus for any $G$-graded algebra  $A$ and a $G$\nobreakdash-module $M$, we obtain a Hopf triple $(A,kG,M)$.
  \end{itemize}
  \end{example}
  
   Given a left Hopf triple $(A, \mathcal{H},M)$, let $C_n(A,M)=M\ot A^{\ot(n+1)}$. 
   We define the simplicial and cyclic operators  on $\{C_n(A,M)\}_n$ by 
  \begin{align}
 &\delta_0(m\ot a_0\ot a_1\ot\dots\ot a_n)=m\ot a_0a_1\ot a_2\ot\dots \ot a_n,\notag\\
&\delta_i(m\ot a_0\ot a_1\ot\dots\ot a_n)=m\ot a_0\ot \dots \ot a_i a_{i+1}\ot\dots \ot  a_n, &&1\le i\le n-1,\notag\\
& \delta_n(m\ot a_0\ot a_1\ot\dots\ot a_n)=a_n^{(-1)}m\ot a_n^{(0)}a_0\ot a_1\dots \ot a_{n-1},\\        
 &\sigma_i(m\ot a_0\ot a_1\ot\dots\ot a_n)=m\ot a_0 \ot \dots \ot a_i\ot 1 \ot \dots \ot a_n, &&0\le i\le n,\notag\\  
 &\tau(m\ot a_0\ot a_1\ot\dots\ot a_n)=a_n^{(-1)}m\ot a_n^{(0)}\ot a_0\ot\dots \ot a_{n-1}.\notag
 \end{align}
 \begin{proposition}
 Endowed with the above operators,  $\{C_n(A,M)\}_{n}$ is a paracyclic module. 
 \end{proposition}
 
 Next, we define a left  $\mathcal{H}$-coaction $\rho:C_n(A,M)\longrightarrow \mathcal{H}\ot C_n(A,M)$ by 
 $$\rho(m\ot a_0\ot\dots \ot a_n)=(a_0^{(-1)}\dots a_n^{(-1)})\ot  m\ot a_0^{(0)}\ot \dots \ot a_n^{(0)}.$$
\begin{lemma}
Endowed with the above coaction,  $C_n(A,M)$ is an   $\mathcal{H}$\nobreakdash-comodule.
\end{lemma}
To define the space of coinvariants  of $C_n(A,M)$, we fix a grouplike element $\sigma\in\mathcal{H}$.
Let
$$C_n^\mathcal{H}(A,M)= C_n(A,M)^{\text{co}\mathcal{H}}=\{ x\in C_n(A,M)\mid \rho(x)=\sigma\ot x\},$$
be the space of coinvariants  of $C_n(A,M)$ with respect to $\sigma$. We would like to find conditions that 
guarantee $\{C_n^\mathcal{H}(A,M)\}_n$ is a cyclic module. This leads us to the following definitions and results.
 
 \begin{definition}\label{definition3.5}
 Let $M$ be a left $\mathcal{H}$\nobreakdash-module and $\sigma\in \mathcal{H}$ a grouplike element.
  We define the $(M,\sigma)$\nobreakdash-{\it twisted antipode} $\widehat{S}:M\ot \mathcal{H}\longrightarrow M\ot\mathcal{H}$ by
  $$\widehat{S}(m\ot h)=h^{(2)}m\ot \sigma S(h^{(1)}).$$
 \end{definition}
 \begin{definition}\label{def3.4}
 Let $M$ be a left $\mathcal{H}$\nobreakdash-module and $\sigma\in\mathcal{H}$ a grouplike element. We call
  $(M,\sigma)$ a 
 matched pair if $\sigma m=m$ for all $m\in M$. We call the matched pair $(M,\sigma)$ 
 a matched pair in involution if 
  $$(\widehat{S})^2=id: M\ot\mathcal{H}\longrightarrow M\ot \mathcal{H},$$
  where $\widehat{S}$ is defined in Definition {\ref{definition3.5}}.
 \end{definition}
 \begin{example} \label{example3.6}
Let $M=k_\delta$ be  the one dimensional module defined by a character $\delta\in \mathcal{H}$. It is 
 clear that $(M,\sigma)$ is a matched pair in involution if and only if $(\delta, \sigma)$
  is a {\it modular pair in involution } 
  in the sense of  \cite{kr2}, i.e.,  $\delta(\sigma)=1$ and $(\sigma\widetilde{S}_\delta)^2=id$.
  \end{example}
  
  Let $k$ be a field of characteristic zero and $q\in k$, $q\neq 0$ and $q$ 
not a root of unity. The Hopf algebra $\mathcal{H}=A(SL_q(2,k))$ is defined as follows.
 As an algebra it is generated by symbols $x,\; u,\; v,\; y,$ with the following relations:
$$ ux=qxu, \;\; vx=qxv, \;\; yu=quy,\;\;yv=qvy,$$
$$uv=vu,\;\; xy-q^{-1}uv=yx-quv=1. $$
The coproduct, counit and antipode of   $\mathcal{H}$ are defined by  
$$\Delta (x)=x \otimes x+u \otimes v,\;\;\;\Delta (u)=x \otimes u+u \otimes y, $$ 
$$\Delta (v)=v \otimes x+y \otimes v,\;\;\;\Delta (y)=v \otimes u+y \otimes y, $$
$$\epsilon (x)=\epsilon (y)=1,\;\;\;\epsilon (u)=\epsilon (v)=0,$$
$$S(x)=y,\;\; S(y)=x,\;\;S(u)=-qu,\;\;S(v)=-q^{-1}v. $$
  For more details about $\mathcal{H}$ we refer to ~\cite{kl}.
  We give an example of a  Hopf triple where  $M$ is not one dimensional.
\begin{example}\label{example3.7}
Let $M$ be a free $k$\nobreakdash-module generated by $m_1$ and $m_2$ and let $\mathcal{H}=A(SL_q(2,k))$ act
 on $M$ as follows:
\begin{center}
$xm_1=qm_2$,\qquad\qquad$xm_2=qm_1$,\\
$um_1=um_2=vm_1=vm_2=0$,\\
$ym_1=q^{-1}m_2$,\qquad\qquad$ym_2=q^{-1}m_1,$
\end{center}
One can check that $(M,1)$ is a matched pair in involution. 
\end{example}
 The following lemma will play an important role in the proof of  Theorem \ref{1}.
 Its  proof is elementary and hence  omitted.  
 \begin{lemma}
 Let $(M,\sigma)$ be a matched pair. 
 Then $\widehat{S}:M\ot \mathcal{H}\longrightarrow M\ot \mathcal{H}$ is invertible and we have 
 $\widehat{S}^{-1}(m\ot h)=h^{(1)}m\ot S^{-1}(h^{(2)})\sigma.$
\end{lemma}
\begin{definition}
Let $(A,\mathcal{H},M)$ be  a Hopf triple and $\sigma\in \mathcal{H}$ a grouplike element.
 We say $(A,\mathcal{H},M)$ is $\sigma$-compatible  if $(M,\sigma)$ is a matched pair in involution.
\end{definition}

 \begin{lemma}\label{lem3.10}
 Let $(A, \mathcal{H},  M)$ be a $\sigma$-compatible left Hopf triple. Then for any $a\in A$ and $m\in M$ 
 $$a^{(-1)}\sigma S(a^{(-3)})\ot a^{(-2)}m\ot a^{(0)}=\sigma\ot a^{(-1)}m\ot a^{(0)}.$$
 \end{lemma}
 \begin{proof}
 Since $\widehat{S}^{-2}=id$,  we have 
  \begin{multline*}
 a^{(-1)}\sigma S(a^{(-3)})\ot a^{(-2)}m\ot a^{(0)}=\widehat{S}^{-2}(a^{(-1)}\sigma 
 S(a^{(-3)})\ot a^{(-2)}m)\ot a^{(0)}\\
 =\sigma^{-1}S^{-2}(a^{(-1)}\sigma S(a^{(-3)}))\ot S^{-1}(a^{(-1)}\sigma S(a^{(-7)}))
 \sigma a^{(-3)}S(a^{(-5)})a^{(-4)}m\ot a^{(0)}\\
 =\underbrace{\sigma^{-1}S^{-2}(a^{(-2)})\sigma} S^{-1}(a^{(-4)})\sigma\ot a^{(-5)}
 \sigma^{-1}\underbrace{S^{-1}(a^{(-1)})\sigma a^{(-3)}m}\ot a^{(0)}\\
 =a^{(-1)}S^{-1}(a^{(-2)})\sigma\ot a^{(-3)}m\ot a^{(0)}\\
 = \sigma\ot a^{(-1)}m\ot a^{(0)}.
 \end{multline*}
 
 \end{proof}
 The following theorem is  the main result of this section.
 \begin{theorem}\label{1}
 Let $(A, \mathcal{H}, M)$ be a $\sigma$-compatible Hopf triple. 
 Then $\{ C_n^\mathcal{H}(A,M)\}_{n}$  endowed with  simplicial and 
 cyclic operators induced by \text{(1)}, is a 
 cyclic module. 
 \end{theorem}
 \begin{proof}
 As a first step we show that the induced simplicial and cyclic  operators are well
  defined on $\{C_n^\mathcal{H}(A,M)\}_n$ . We just  prove this for  $\tau$, and  $\delta_n$ 
  and leave the rest to the reader. Let $(m\ot a_0\ot\dots \ot a_n)\in C_n^\mathcal{H}(A,M)$. We have 
  \begin{equation}
     a_0^{(-1)}\dots a_n^{(-1)}\ot m\ot a_0^{(0)}\ot \dots a_n^{(0)}=\sigma\ot m\ot a_0\ot \dots a_n 
     \end{equation}
  which implies 
 $$
 a_0^{(-1)}\dots a_{n-1}^{(-1)}\ot m\ot a_0^{(0)}\ot \dots  \ot a_{n-1}^{(0)}\ot a_n=
  \sigma S(a_n^{(-1)})\ot m\ot a_0\ot \dots a_n^{(0)}$$
  and 
  \begin{multline*}
    a_n^{(-1)}a_0^{(-1)}\dots a_{n-1}^{(-1)}\ot a_n^{(-2)}
  m\ot a_n^{(0)}\ot a_0^{(0)}\ot \dots \ot a_{n-1}^{(0)}=\\
  \shoveright{~~~~=\underbrace{a_n^{(-1)}\sigma S(a_n^{(-3)})\ot a_n^{(-2)}m\ot 
  a_n^{(0)}}\ot a_0\ot \dots \ot a_{n-1}.}\\
  \end{multline*}
  Applying   Lemma \ref{lem3.10} for $a=a_n$ we have 
  \begin{multline}
    a_n^{(-1)}a_0^{(-1)}\dots a_{n-1}^{(-1)}\ot a_n^{(-2)}
  m\ot a_n^{(0)}\ot a_0^{(0)}\ot \dots \ot a_{n-1}^{(0)}=\\
  \sigma \ot a_n^{(-1)}m\ot a_n^{(0)}\ot a_0\ot \dots \ot a_{n-1}
  \end{multline}
  which means  $\tau(m\ot a_0\ot\dots \ot a_n)\in C_n^\mathcal{H}(A,M)$.\\ 
   
    From ($3$) we obtain  
    \begin{multline*}
    a_n^{(-1)}a_0^{(-1)}\dots a_{n-1}^{(-1)}\ot a_n^{(-2)}
  m\ot a_n^{(0)} a_0^{(0)}\ot \dots \ot a_{n-1}^{(0)}=\\
  \sigma \ot a_n^{(-1)}m\ot a_n^{(0)} a_0\ot \dots \ot a_{n-1}
  \end{multline*}
   which implies  $d_n(m\ot a_0\ot\dots \ot a_n)\in C_{n-1}^\mathcal{H}(A,M)$.\\
   Checking  that the  other simplicial operators  are well defined on $\{C_n^\mathcal{H}(A,M)\}_n$  is
    straightforward. 
    
   The only thing left is to show that $\tau^{n+1}=id$. We have    
   $$\tau^{n+1}(m\ot a_0\ot\dots \ot a_n)=a_0^{(-1)}\dots a_n^{(-1)} m\ot a_0^{(0)}\ot \dots a_n^{(0)}.$$   
   Now since we are in $C_n^\mathcal{H}(A,M)$,  and by ($2$) we have  
   $$\tau^{n+1}(m\ot a_0\ot\dots \ot a_n)=\sigma m\ot a_0\ot\dots a_n=m\ot a_0\ot\dots \ot a_n,$$
    because $(\sigma, M)$ is a matched pair.
   \end{proof}
   
   We denote the resulting Hochschild, cyclic and periodic cyclic homology groups 
   of the cyclic module $\{C_n^\mathcal{H}(A,M)\}_n$ by $HH_\bullet^\mathcal{H}(A,M)$,
    $HC_\bullet^\mathcal{H}(A,M)$ and
   $HP^\mathcal{H}_\bullet(A,M)$, respectively, and refer to them  as {\it invariant
    Hochschild}, {\it cyclic} and 
   {\it periodic cyclic } homology groups of the $\sigma$-compatible Hopf triple $(A,\mathcal{H},M)$.
   
   We give a few  examples of invariant cyclic homology. 
   More examples can be  found in Subsections \ref{subsection3.2},  \ref{subsection3.3}, and  \ref{subsection3.4}. 
   It is clear that if $(A,k,k)$ is a trivial Hopf triple 
   ( Example \ref{example3.2} (i)), then $HC^k_\bullet(A,k)\cong HC_\bullet(A)$,
    i.e.,  in this case,  invariant cyclic homology is  the same  as cyclic homology of algebras.
    
      Next we  consider the Hopf triple $(A,kG,k)$ defined 
	in Example \ref{example3.2}(iv). Computing the invariant cyclic homology group  $HC_\bullet^{kG}(A,k)$, except in
	especial cases, is not easy in general. Let $G=\mathbb{Z}$ and $A=\underset{n\ge 0}{\oplus}A_n$
	 be a positively graded algebra. One can see that 
	 $$C_n^{kG}(A,k)\cong C_n(A_0)\;\;n\ge 0.$$
	 Therefore, we obtain 
	 $$HC_n^{k\mathbb{Z}}(A,k)\cong HC_n(A_0).$$
	  This statement is false if $A$ has nonzero negative components. 
	  For example, if $A=k\lbrack z,z^{-1}\rbrack =k\mathbb{Z}$ is the (Hopf) algebra of Laurent polynomials,
	  by Proposition \ref{prop3.14}, we have 
	  $HC_n^{k\mathbb{Z}}(k\mathbb{Z},k)\cong \widetilde{HC}_n^{(\epsilon,1)}(k\mathbb{Z})=k$  for    
	  $n\ge 0.$
    
    The following lemma  enables us   to identify 
    $C_n^\mathcal{H}(\mathcal{H},M)$ with $M\ot \mathcal{H}^{\ot{n}}.$
    \begin{lemma}\label{lemma3.12}
    Let $V$ be a left $\mathcal{H}$\nobreakdash-comodule. Then $(\mathcal{H}\ot V)^{\text{co}\mathcal{H}}\cong V$.
    \end{lemma}
    \begin{proof}
    Define $\eta: V\longrightarrow ( \mathcal{H}\ot V)^{\text{co}\mathcal{H}}$  by
     \begin{center}
     { $\eta(v)=\sigma S(v^{(-1)})\ot v^{(0)}.$ }
\end{center}
 We show that   $\eta$ is well defined and is an isomorphism of $k$\nobreakdash-modules.
  Let $\rho$ be the structure map of $V$ and $\bar\rho$ be the induced comodule structure on $\mathcal{H}\ot V$.
   Since 
 \begin{center}
 $\bar\rho(\eta(v))=\bar\rho(\sigma S(v^{(-1)})\ot v^{(0)})=\sigma S(v^{(-2)})v^{(-1)}\ot\sigma S(v^{(-3)})\ot v^{(0)}
 $\\$=\sigma\ot\sigma S(v^{(-1)})\ot v^{(0)},$
 \end{center}
  we see that $\eta$ is  well defined. Now consider the following map: 
  
 \begin{center}
{    $\theta:(\mathcal{H}\ot V)^{\text{co}\mathcal{H}} \longrightarrow V$ }\\
{ $\theta(h\ot v)=\epsilon(h)v.$ }
\end{center} 
It is easy to see that $\theta\eta=id_V$. We complete the proof of this lemma by showing 
that $\eta\theta=id_{\mathcal{H}\ot V}$. Let $h\ot v\in (\mathcal{H}\ot V)^{\text{co}\mathcal{H}}.$ Then 
\begin{multline*}
\eta(\theta(h\ot v))=\eta(\epsilon(h)v)=\sigma S(v^{(-1)})
\epsilon(h)\ot v^{(0)}=\\=\sigma S(h^{(1)}v^{(-1)})h^{(2)}\ot v^{(0)}=\sigma S(\sigma)h\ot v=h\ot v.
\end{multline*}
    \end{proof}
Applying the above isomorphism to $V=M\ot \mathcal{H}^{\ot(n+1)}$, one can identify 
$C_n^{\mathcal{H}}(\mathcal{H},M)$ with 
 $M\ot \mathcal{H}^{\ot{n}}$ and its simplicial and cyclic operators as follows:
 \begin{align}{\label{hasan}}
&\delta_0(m\ot h_1\ot \dots\ot h_n)=\epsilon(h_1)m\ot h_2\ot\dots\ot h_n,\notag\\
&\delta_i(m\ot h_1\ot \dots\ot h_n)=m\ot h_1\ot\dots \ot h_ih_{i+1}\ot\dots\ot 
h_n,&& \hspace{-34pt}1\le i\le n-1\notag\\
&\delta_n(m\ot h_1\ot \dots\ot h_n)=h_nm\ot h_1\ot\dots\ot h_{n-1},\notag\\
&\sigma_0(m\ot h_1\ot \dots\ot h_n)=m\ot 1\ot h_1\ot\dots\ot h_n,\\
&\sigma_i(m\ot h_1\ot \dots\ot h_n)=m\ot  h_1\ot\dots \ot h_i\ot 1\ot \dots\ot h_n, &&\hspace{-34pt}{1\le i\le n}\notag\\
&\tau(m\ot h_1\ot \dots\ot h_n)=h_n^{(2)}m\ot \sigma S(h_1^{(1)}\dots h_n^{(1)})\ot h_1^{(2)}\ot\dots\ot
 h_{n-1}^{(2)}.\notag
\end{align}

We see that the Hochschild complex of $\{C_n^\mathcal{H}(\mathcal{H},M)\}_n$ is 
isomorphic to the Hochschild complex of the algebra $\mathcal{H}$ with
 coefficients in the bimodule $M$ where the right action is via $\epsilon$.
 
Let us  recall that a Hopf algebra is called {\it semisimple} 
    if the underlying algebra of $\mathcal{H}$ is semisimple  \cite{kl}. 
  It is known that  $\mathcal{H}$ is semisimple if and only if there is
   a {\it normalized integral} in $\mathcal{H}$, i.e., an element 
  $t\in \mathcal{H}$ such that $th=\epsilon(h)t$ for all $h\in \mathcal{H}$ and $\epsilon(t)=1 $.
\begin{proposition}
   Let $(\mathcal{H},\mathcal{H},M)$ be a $\sigma$-compatible Hopf triple. If $\mathcal{H}$ is semisimple then 
$HC_{2n}^{\mathcal{H}}(\mathcal{H},M)=M_{\mathcal{H}}$ and 
$HC_{2n+1}^{\mathcal{H}}(\mathcal{H},M)=0$ for all $n\ge 0$.
\end{proposition}
\begin{proof}
Using  a normalized integral, we  define the following homotopy operator
      $h:M\ot\mathcal{H}^{\ot n} \rightarrow M\ot\mathcal{H}^{\ot n+1}$ by 
  $$h(m\ot h_1\ot \dots \ot h_n)=m\ot t\ot h_1\ot \dots \ot h_n.$$  
 It can be checked that  $hb+bh=id$, where $b=\sum_{i=0}^n(-1)^i\delta_i $ is the Hochschild boundary map.
  It follows that $H_i^\mathcal{H}(\mathcal{H},M)=0$ for $i\ge 1$, and $H_0^\mathcal{H}
  (\mathcal{H},M)=M_\mathcal{H}$. The rest is obvious.
  \end{proof}
 
 As another example we like to compute  the invariant cyclic homology  of the  
 Hopf triple introduced in Example {\ref{example3.7}}. We  first compute 
 the Hochschild homology $H_\bullet(\mathcal{H}, M)$. 
One knows that $H_{\bullet}(\mathcal{H}, M)=Tor_{\bullet}^{\mathcal{H}^e}(\mathcal{H}, M)$, where $\mathcal{H}^e =
 \mathcal{H}\otimes \mathcal{H}^{op}$.  We take advantage of the  free resolution for
 $\mathcal{H}$ given in ~\cite{ma}: 
$$\dots \rightarrow M_2 \overset{d_2}{\rightarrow} M_1 \overset{d_1}{\rightarrow} 
M_0 \overset{\mu}{\rightarrow} \mathcal{H},$$   
  where $\mu$ is the augmentation map and $M_\ast$,  is a family of free left $\mathcal{H}^e$-modules  
  with their  rank  given by 
  \begin{eqnarray*}
 && rank(M_0)=1\\
 && rank(M_1)=4\\
 && rank(M_2)=7\\
 && rank(M_\ast)=8,\; \ast\ge 3.
   \end{eqnarray*}
 We refer the interested reader to \cite{ma} for this resolution.
  After a lengthy computation we obtain  the following 
 theorem.
 \begin{theorem}
For any $q\in k$ which is not a root  of unity one has\\
 ${HC}^\mathcal{H}_1(A(SL_q(2)),M)=k\oplus k$ and
  ${HC}^\mathcal{H}_n(A(SL_q(2)),M)= 0$  for all $n\neq 1$.     
\end{theorem}

	  In ~\cite{kr2}, and independently ~\cite{tr}, for  a given Hopf algebra $\mathcal{H}$ endowed with a modular
      pair in involution $(\delta, \sigma)$ in the sense of {\cite{kr2}}, a cyclic module $\{
       \widetilde{H}_n^{(\delta,\sigma)}\}_{n}$
   is defined. The resulting cyclic homology theory is, in a sense, dual to the cyclic cohomology 
   theory of Hopf algebras defined by Connes and Moscovici. We  show  that this theory is an 
   example of invariant cyclic homology theory defined in this section. Consider 
   the $\sigma$-compatible Hopf triple 
   $(\mathcal{H},\mathcal{H},k_\delta)$ defined in Example \ref{example3.2}(iv). One can check  that for $M=k_\delta$ 
   the operators in (\ref{hasan}) are exactly the operators defined in \cite{kr2}. This proves the following
    proposition. 
   \begin{proposition}\label{prop3.14}
   The cyclic modules $\{\widetilde{\mathcal{H}}_n^{(\delta,\sigma)}\}_{n}$
    and $\{C_n^\mathcal{H}(\mathcal{H},k_\delta)\}_n$ are isomorphic. 
   \end{proposition} 
          
We show that under suitable conditions, the invariant cyclic homology is a direct summand in 
cyclic homology of algebras. Consider a $\sigma$-compatible Hopf triple $(A,\mathcal{H},k)$
 where the action of $\mathcal{H}$ on $k$ is via the counit $\epsilon$. It is  easy
to see that in this case the paracyclic module $\{C_n(A,k)\}_n$ is isomorphic to $\{C_n(A)\}_n$, 
the cyclic module of the algebra 
$A$. We therefore obtain an inclusion 
$$i:C_n^{\mathcal{H}}(A,k)\hookrightarrow C_n(A).$$
To define a left inverse for $i$, we need a suitable linear functional on $\mathcal{H}$.
 Recall from \cite{kr2} that a trace 
$Tr: \mathcal{H} \longrightarrow k$ is called $\sigma$-{\it invariant} if for any $h\in \mathcal{H}$
$$Tr(h^{(1)}){h^{(2)}}=Tr(h)\sigma.$$
Given such a trace on $\mathcal{H}$, we define the {\it averaging operator}
\begin{center}
$\gamma: C_n(A)\longrightarrow C_n^{\mathcal{H}}(A)$
$\gamma(a_0\ot \dots\ot a_n)=Tr(a_0^{(-1)}\dots a_n^{(-1)})a_0^{(0)}\ot a_1^{(0)}\ot \dots \ot a_n^{(0)}.$
\end{center}
Let $\rho: A\longrightarrow \mathcal{H}\ot A$ denote the coaction and $\bar\rho$ the induced coaction. We have 
\begin{multline*}
\bar\rho\gamma(a_0\ot \dots \ot a_n)=Tr(a_o^{(-2)}\dots a_n^{(-2)})a_0^{(-1)} a_1^{(-1)}
 \dots  a_n^{(-1)}\ot a_0^{(0)}\ot \dots \ot a_n^{(0)}\\
 =\sigma\ot Tr(a_o^{(-1)})\dots 
 a_n^{(-1)})a_0^{(0)}\ot a_1^{(0)}\ot \dots \ot a_n^{(0)},
 \end{multline*}
which shows that  the image of $\gamma$ is   in fact in the subspace of invariant chains. 
 The following proposition has an elementary proof. 
\begin{proposition}\label{p1}
Let $(A,\mathcal{H},k)$ be  a $\sigma$-compatible Hopf triple and  $\mathcal{H}$ 
admits a $\sigma$-invariant trace $Tr$. Then $\gamma$  is a 
cyclic map and  $\gamma i=Tr(\sigma)id$.
\end{proposition}
\begin{corollary}
If $Tr(\sigma)$ is invertible  in $k$, then $HC_\bullet^{\mathcal{H}}(A,k)$ 
is a direct summand in $HC_\bullet(A)$.
\end{corollary}
Let $(\mathcal{H},\mathcal{H},M)$ be the    Hopf triple defined in Example \ref{example3.2}(ii). 
 One can easily see that  $(M_n(\mathcal{H}),\mathcal{H},M)$ is also a  Hopf triple, where
 the coaction of $\mathcal{H}$ on $M_n(\mathcal{H})$ is induced by the comultiplication of $\mathcal{H}$, i.e., 
 for all $h\ot u\in \mathcal{H}\ot M_n(k)=M_n(\mathcal{H})$,  
{$\rho(h\ot u)=h^{(1)}\ot h^{(2)}\ot u$}. We have the following identifications  
$$C_n(M_k(\mathcal{H}),M)=M\ot M_k(\mathcal{H})^{\ot{(n+1)}}=\mathcal{H}\ot M_k(M) \ot  M_k(\mathcal{H})^{\ot{n}}.$$
  Now  we can apply Lemma \ref{lemma3.12} to \mbox{$V= M_k(M)\ot M_k(\mathcal{H})^{\ot{n}}$}, to get 
  $$C_n^\mathcal{H}(M_k(\mathcal{H}),M)\cong M_k(M)\ot M_k(\mathcal{H})^{\ot{n}}.$$
  \begin{theorem}{(Morita invariance)}
  For any matched pair in involution $(M,\sigma)$,   and any  $k\ge 1$ one has 
  $$HC^\mathcal{H}_n(\mathcal{H},M)\cong HC_n^\mathcal{H}(M_k(\mathcal{H}),M),\qquad n\ge 0.$$
    \end{theorem}
\begin{proof}
 Consider the inclusion map 
 $$C_\bullet^\mathcal{H}(\mathcal{H},M)\hookrightarrow C^\mathcal{H}_\bullet(M_k(\mathcal{H}),M). $$
 One can show that this map  is well-defined and is a 
 cyclic module  map. On the other hand by the above explanation 
   one can see that 
 $$H_n^\mathcal{H}(M_k(\mathcal{H}),M)\cong H_n(M_k(\mathcal{H}),M_k(M))\cong
  H_n(\mathcal{H},M)\cong H_n^\mathcal{H}(\mathcal{H},M),$$
  where  the second isomorphism  is Morita invariant of ordinary Hochschild homology  \cite{ld}. The theorem 
  can be proved now by invoking the long exact sequence relating the 
  (invariant) Hochschild and cyclic homology groups. 
\end{proof}
\subsection{Relation with invariant de Rham cohomology}\label{subsection3.2}
In \cite{ce}, Chevalley and Eilenberg defined the invariant de Rham  cohomology of a $G$-manifold
where $G$ is a Lie group. One is naturally interested to know to what extent  the relation 
between Hochschild  homology  and differential forms and between cyclic homology and de Rham  cohomology 
 (Hochschild-Kostant-Rosenberg \cite{hkr} and Connes \cite{ac85}) 
 extend to our invariant setting. While we do not have a proof, we believe
   these   results can be extended to  the invariant case and in this subsection 
   present some evidence in this direction.
   
 Let $\mu: G\times V\longrightarrow V$ denote the action of an  affine 
 algebraic  group on a smooth affine algebraic variety  $V$.
 Let $\mathcal{H}=\mathbb{C}\lbrack G\rbrack$ and $A=\mathbb{C}\lbrack V\rbrack$  
 be the coordinate rings of 
$G$ and $V$, respectively. Then $A$ is an $\mathcal{H}$\nobreakdash-comodule algebra via a map
$\rho:A\longrightarrow \mathcal{H}\ot A$ which  is obtained by dualizing  $\mu$. Let 
$\Omega^\bullet A$  denote the algebraic de Rham complex of $V$ and
$\Omega^\bullet_{\text{inv}} A=\{\omega\in \Omega^\bullet A\mid g^\ast\omega=\omega,\;\forall g\in G \}$
its invariant part.  The {\it invariant de Rham cohomology  }of $V$ \cite{ce}
is, by definition, the cohomology of the complex $\Omega_{\text{inv}}^\bullet A$. Consider 
the map $\pi: A^{\ot{(n+1)}}\longrightarrow \Omega^n A$, 
\begin{center}
$\pi(a_0\ot a_1\ot\dots \ot a_n)=a_0da_1\dots da_n$,
\end{center} 
and the antisymmetrization map $\gamma: \Omega^nA\longrightarrow HH_n(A)$,
\begin{center}
$\gamma(a_0da_1\dots da_n)=\frac{1}{n!}
\sum_{\alpha\in S_n}sgn(\alpha)(a_0\ot a_{\alpha(1)}\ot\dots\ot a_{\alpha(n)}).$
\end{center}
It is easy to check that $\pi$ and $\gamma$ induce maps \\
\begin{eqnarray*}
&&\pi^G:(A^{\ot(n+1)})^{\text{co}\mathcal{H}}\longrightarrow \Omega_{\text{inv}}^n A\\
&&\Omega_{\text{inv}}^n A\longrightarrow HH_n^{\mathcal{H}}(A,k).
\end{eqnarray*}
 We therefore obtain maps
 $$HH_n^\mathcal{H}(A,k)\longrightarrow \Omega^n_{\text{inv}}(A)$$
 { and}
 $$HP^\mathcal{H}_n(A,k)\longrightarrow \underset{i=n \;{\text{mod}}\;2}
 {\bigoplus}H_{\text{dR, inv}}^i(V).$$
 
 We believe that  both maps are isomorphisms, but do not have a proof at this stage except for $V=G$
  with translation action. In this case the Hochschild homology of the 
  cyclic module $\{C_n^\mathcal{H}(\mathcal{H},\mathbb{C})\}_n$ is isomorphic to $H_n(\mathcal{H},\mathbb{C})$.
   By  Hochschild-Kostant-Rosenberg's theorem \cite{hkr}, we obtain
   $$H_n^\mathcal{H}(\mathcal{H},\mathbb{C})\cong \land^n(\mathfrak{g}^\ast),$$  where $\mathfrak{g}=\text{Lie}(G)$
    is the Lie algebra
   of $G$.
One can then check that 
$$HP_n^\mathcal{H}(\mathcal{H},\mathbb{C})\cong
\bigoplus_{i=n \;{\text{mod}}\;2}H^i_{\text{Lie}}(\mathfrak{g},\mathbb{C})\cong
\bigoplus_{i=n \;{\text{mod}}\;2}H^i_{\text{dR, inv}}(G).$$

 \subsection{The cocommutative case}\label{subsection3.3}
In this subsection we  show that Lie algebra 
homology (with  coefficients)  and group homology (with  coefficients) can be interpreted in  
our framework. In fact we prove a much more general result  for all cocommutative Hopf algebras. 
Consider the Hopf triple
 $(\mathcal{H},\mathcal{H},M)$  defined in Example \ref{example3.2}(ii).  We will show that if $\mathcal{H}$
  is cocommutative, then
 $$HC_n^{\mathcal{H}}(\mathcal{H},M)=\bigoplus_{i\ge 0} H_{n-2i}(\mathcal{H},M)$$ 
$$HP_n^{\mathcal{H}}(\mathcal{H},M)=\bigoplus_{i=n\;\text{mod}\;2} H_i(\mathcal{H},M)$$ 
 where in the right hand side we have {\it Hopf homology} groups. For $\mathcal{H}=U(\mathfrak{g})$
  or $\mathcal{H}=kG$ we obtain the relation 
between invariant cyclic homology on one side and Lie algebra or group homology on the other side.

We recall the definition of   Hopf homology groups .
  This  notion extends group homology  and  Lie algebra homology  
  with coefficients in a module to all Hopf 
algebras. 
 Let $\mathcal{H}$ be a Hopf algebra.  The ground ring $k$ is a right $\mathcal{H}$\nobreakdash-module via 
 the map $(r,h)\mapsto r\epsilon(h)$, for all $ r\in k$, and $h\in\mathcal{H}$. 
 It is clear that the functor of coinvariants $M\mapsto M_\mathcal{H}=k \underset{\mathcal{H}}{\ot}M$ 
 from the category of left $\mathcal{H}$\nobreakdash-modules to the category of $k$\nobreakdash-modules is right
  exact.
 We denote the corresponding left derived functors by $H_n(\mathcal{H},M))$.
 For $\mathcal{H}=U(\mathfrak{g})$ or $\mathcal{H}=kG$ we obtain Lie algebra or group homology, 
 respectively.
 
 %We first recall the homology of a Hopf algebra with coefficients in a module. Let $\mathcal{H}$ be a Hopf algebra 
 %and $M$ a left $\mathcal{H}$\nobreakdash-module, then $M$ is an $\mathcal{H}$-bimodule, where the right action is 
%via 
 %the counit $\epsilon$. We  have then the Hochschild homology group $H_\bullet(\mathcal{H},M)$. Alternatively,
  %these groups are higher derived functions of the functor of invariants 
  %$\mathcal{H}-\text{mod}\longrightarrow k-\text{mod}$, $M\mapsto m^\mathcal{H}:=\{m\in M\mid hm=\epsilon(h)m\}.$ 
   
    Recall that the {\it path space }    $EM$ of a simplicial module $M=\{M_n\}_{n}$
	 is defined by $(EM)_n=M_{n+1}$. It is a simplicial module with the simplicial operators of 
	 $M$ shifted by one. If $M$ is a cyclic module, there is no natural cyclic structure on $EM$. 
	 It is therefore remarkable that if $\mathcal{H}$ is cocommutative the path space
	   $\{EC^\mathcal{H}(\mathcal{H},M)\}_n$ is a  cyclic module in a natural way. Define the operators
	\begin{center}
	$t_n:EC^\mathcal{H}_n(\mathcal{H},M)\longrightarrow EC^\mathcal{H}_n(\mathcal{H},M), \qquad $ \\
  \makebox[8in][l]{by}
	$t_n(m\ot h_0\ot \dots \ot h_n )=h_n^{(3)}m\ot h_0h_1^{(1)}\dots 
	h_n^{(1)}\ot S(h_1^{(2)}\dots h_n^{(2)})\ot h_1^{(3)}\ot \dots \ot h_{n-1}^{(3)}.$
	\end{center}
	The proof of the following proposition is similar to Lemma $4.1$ in \cite{kr2} hence we discard it.
	 \begin{proposition}
	 Let $\mathcal{H}$ be a cocommutative Hopf algebra. Then the \mbox{  path space}
	  $\{EC_n^\mathcal{H}(\mathcal{H},M)\}_n$  is a cyclic module.
	 \end{proposition}   
	 We define a map 
	 \begin{center}
	 $\theta: EC_n^\mathcal{H}(\mathcal{H},M)\longrightarrow C_n^\mathcal{H}(\mathcal{H},M) $\\
	 $\theta(m\ot h_0\ot \dots \ot h_n )=\epsilon(h_0)m\ot h_1\ot \dots \ot h_n. $
	 \end{center}
	  In a natural way $\mathcal{H}$ 
	 has an  action on $EC_n^\mathcal{H}(\mathcal{H},M)$  defined by 
	 $g(m\ot h_0\ot \dots \ot h_n )=m\ot gh_0\ot \dots \ot h_n $.
	 \begin{lemma}If $\mathcal{H}$ is a cocommutative Hopf algebra then
	  $\{EC_n^\mathcal{H}(\mathcal{H},M)\}_n$  is a cyclic module and $\theta $ is a cyclic module map. Moreover 
    $$k\underset{\mathcal{H}}{\ot}EC_n^\mathcal{H}(\mathcal{H},M)\cong C_n^\mathcal{H}(\mathcal{H},M).$$
	 \end{lemma} 
	 Using the above lemma and the method used in the proof of Theorem $4.1$ in \cite{kr2}, we can prove 
	 the following theorem.  
	 
   \begin{theorem}\label{hasan2}
   Let $\mathcal{H}$ be a cocommutative Hopf algebra. 
   Then $$HC_n^\mathcal{H}(\mathcal{H}, M)\cong \underset{i\ge 0}{\bigoplus}H_{n-2i}(\mathcal{H},{M}),$$
    where in the  right
    hand side appears 
    the Hopf homology groups. 
   \end{theorem}

   \begin{corollary}
    Let $\mathfrak{g}$ be a Lie algebra and $M$ a $\mathfrak{g}$\nobreakdash-module. Then\\
      \mbox{$HC_n^{U(\mathfrak{g})}(U(\mathfrak{g}),M)\cong \underset{i\ge
      0}{\bigoplus}H_{n-2i}(\mathfrak{g},M)$,}
      where in the right  hand we have the Lie algebra homology. 
   \end{corollary}
   %\begin{proof}
   %By Theorem, $HC_n^{U(\mathfrak{g})}(U(\mathfrak{g}),M)\cong \underset{i\ge 0}
   %{\bigoplus}H_{n-2i}(U(\mathfrak{g});{M})$. 
   %Now, it well known that $H(U(\mathfrak{g});\widetilde{M})\cong H_n(\mathfrak{g},M)$ ~\cite{ld}.
   %\end{proof}
  %By a similar argument, we obtain
  \begin{corollary}
  Let $G$ be a (discrete) group and $M$ a $G$\nobreakdash-module. Then 
  $HC_n^{kG}(kG,M)\cong \underset{i\ge 0}{\bigoplus}H_{n-2i}(G,M)$, where in the right hand side we have   group
   homology.
  \end{corollary}            
        
\subsection{Invariant cyclic homology of smash products}\label{subsection3.4}
Let $B$ be a right $\mathcal{H}$\nobreakdash-module algebra and let $A=\mathcal{H}\# B$ be the smash product
of $\mathcal{H}$ and $B$. Recall that, as a  $k$\nobreakdash-module, $A=\mathcal{H}\ot B$, and the
product of $A$ is defined by $$(h\ot a )(g\ot b)=hg^{(1)}\ot (a)g^{(2)}b.$$
Now $A$ is a left  $\mathcal{H}$\nobreakdash-comodule algebra with coaction $\rho:A\longrightarrow  \mathcal{H}\ot A$ 
defined by 
           $$\rho(h\ot a )=h^{(1)}\ot (h^{(2)}\ot a).$$ 
Let $M$ be a left $\mathcal{H}$\nobreakdash-module and $\sigma\in \mathcal{H}$
a grouplike element such that $(A,\mathcal{H},M)$                
  is a $\sigma$-compatible Hopf triple. In this subsection we give a spectral 
  sequence to compute
  the invariant cyclic homology groups $\{HC_n^\mathcal{H}(A,M)\}_n$. 
  
  The main idea  is  to define a cylindrical module 
  $\{X_{p,q}\}_{p,q}$, in such a way that its diagonal be  isomorphic 
   to the cyclic module 
   $\{C^\mathcal{H}_n(A,M)\}_n$ (cf. \cite{ak} for a similar spectral sequence 
   for Hopf algebra equivariant cyclic homology which motivated this subsection). We can then use the Eilenberg-Zilber
    theorem for cylindrical modules to obtain the desired spectral sequence.  Let 
     $$X_{p,q}=M \ot \mathcal{H}^{\ot p}\ot B^{\ot (q+1)}.$$  We define the
      vertical and horizontal 
      simplicial and cyclic operators by\\
   $\delta_0(m\ot g_1\ot \dots \ot g_p \ot a_0\ot\dots \ot a_q)=m\ot g_1\ot \dots \ot g_p
    \ot a_0a_1\ot a_2\ot \dots \ot a_q$\\
    \\
   $\delta_i(m\ot g_1\ot \dots \ot g_p \ot a_0\ot\dots \ot a_q)=m\ot g_1\ot \dots \ot g_p \ot a_0\ot \dots 
   \ot a_ia_{i+1}\ot \dots \ot a_q$\\
   \\
   $\delta_q(m\ot g_1\ot \dots \ot g_p \ot  a_0\ot\dots \ot a_q)=m\ot g_1^{(2)}\ot \dots \ot g_p^{(2)}
    \ot$\\
    $ ~~~~~~~~~~~~~~~~~~~~~~~~~~~~~~~~~~~~~~
    (a_q)\sigma S(g_1^{(1)} \dots g_p^{(1)})g_1^{(3)} \dots g_p^{(3)} a_0\ot a_1\ot  \dots \ot a_{q-1}$\\
    \\
    $\tau (m\ot g_1\ot \dots \ot g_p \ot  a_0\ot\dots \ot a_q)=m\ot g_1^{(2)}\ot \dots \ot g_p^{(2)}
    \ot$\\
    $ ~~~~~~~~~~~~~~~~~~~~~~~~~~~~~~~~~~~~~~
    (a_q)\sigma S(g_1^{(1)} \dots g_p^{(1)})g_1^{(3)} \dots g_p^{(3)}\ot a_0\ot a_1\ot  \dots \ot a_{q-1}$\\
    \\
    $d_0(m\ot g_1\ot \dots \ot g_p \ot a_0\ot\dots \ot a_q)=\epsilon(g_1)m\ot g_2\ot \dots \ot g_p
    \ot a_0\ot a_1\ot  \dots \ot a_q$\\
    \\
    $d_i(m\ot g_1\ot \dots \ot g_p \ot a_0\ot\dots \ot a_q)=m\ot g_1\ot\dots \ot g_ig_{i+1}\ot \dots \ot g_p
    \ot a_0\ot a_1\ot  \dots \ot a_q$\\
    \\
    $d_p(m\ot g_1\ot \dots \ot g_p \ot a_0\ot\dots \ot a_q)=g_p^{(1)}m\ot g_1\ot \dots \ot g_{p-1}
    \ot$\\
    $~~~~~~~~~~~~~~~~~~~~~~~~~~~~~~~~ (a_0)S^{-1}(g_p^{(q+2)})\ot  \dots \ot (a_q)S^{-1}(g_p^{(2)})$\\
    \\
    $t(m\ot g_1\ot \dots \ot g_p \ot a_0\ot\dots \ot a_q)=g_p^{(2)}m\ot 
    \sigma S(g_1^{(1)} \dots g_{p}^{(1)})\ot g_1^{(2)}\ot \dots \ot g_{p-1}^{(2)}
    \ot$\\
    $~~~~~~~~~~~~~~~~~~~~~~~~~~~~~~~~ (a_0)S^{-1}(g_p^{(q+3)})\ot  \dots \ot (a_q)S^{-1}(g_p^{(3)}).$
    
    The proof of the following theorem involves several pages of verifications and we simply omit it.
   \begin{theorem}
Endowed with the operators defined above,  $\{X_{p,q}\}_{p,q}$ is a cylindrical module. 
\end{theorem}
 
 The diagonal of this cylindrical module can be identified with the invariant cyclic module
  $\{C_n^\mathcal{H}(A,M)\}_n$:
 \begin{proposition}\label{proposition3.25}
 There exists a natural isomorphism of cyclic modules $d(X)\cong \{C_n^\mathcal{H}(A,M)\}_n$.
 \end{proposition}  
 \begin{proof}
 We have
 \begin{align*}
  C_n^\mathcal{H}(A,M)=&(M\ot A^{\ot(n+1)})^{\text{co}\mathcal{H}}\\
  &=\lbrack M\ot
  (\mathcal{H}\ot B)^{\ot(n+1)}\rbrack^{\text{co}\mathcal{H}}\cong M\ot B\ot (\mathcal{H}\ot B)^{\ot n},
  \end{align*}
   where in the last isomorphism we have used Lemma \ref{lemma3.12}.
 One can check that the following two maps are cyclic maps  and inverse to one another. 
  { $$\phi: M\ot B\ot (\mathcal{H}\ot B)^{\ot (n+1) }\longrightarrow M\ot \mathcal{H}^{\ot n}\ot B^{\ot (n+1)}$$}
  \begin{center}
  $\phi(m\ot a_0\ot g_1\ot a_1\ot\dots \ot g_n\ot a_n)=m\ot g_1^{(1)}\ot \dots
   \ot g_n^{(1)}\ot (a_0)g_1^{(2)}\dots g_n^{(2)}\ot (a_1)g_2^{(3)}\dots g_n^{(3)}
   \ot\dots \ot (a_{n-1})g_n^{(n+2)}\ot a_n$
   $$\psi:M\ot \mathcal{H}^{\ot (n)}\ot B^{\ot (n+1)}\longrightarrow M\ot B\ot (\mathcal{H}\ot B)^{\ot(n+1) }$$
   $\psi(m\ot g_1\ot\dots\ot g_n\ot a_0\ot\dots\ot a_n)=m\ot (a_0)S^{-1}(g_1^{(2)}\dots g_n^{(n+1)})\ot g_1^{(1)}
   \ot (a_1)S^{-1}(g_2^{(2)}\dots g_n^{(n)})\ot \dots \ot g_{n-1}^{(1)}\ot (a_{n-1})
   S^{-1}(g_n^{(2)})\ot g_n^{(1)}\ot a_n.$
   \end{center}
 \end{proof}
 Now we are in a position to apply the Eilenberg-Zilber theorem for cylindrical modules combined with Proposition 
 \ref{proposition3.25} to conclude that
$$ Tot(X)\cong d(X)\cong \{C_n^\mathcal{H}(A,M)\}_n. $$ 
 We  use the following filtration   to derive a spectral sequence for  $HC_\bullet^\mathcal{H}(A,M)$
 $$F_iTotX=\sum_{\underset{q\le i}{p+q=n}}X_{p,q}=\sum_
 {\underset{q\le i}{p+q=n}}M\ot\mathcal{H}^{\ot p}\ot B^{\ot (q+1)}. $$
 \begin{lemma}
 The $E^1$-term of the spectral sequence associated to the above filtration   is given by the Hopf homology
  groups 
 $$E^1_{p,q}=H_p(\mathcal{H},M\ot B^{\ot(q+1)}),$$ where the action of $\mathcal{H}$ on $M\ot B^{\ot(q+1)}$ is
  given by
 $$h(m\ot b_0\ot\dots\ot b_n)=h^{(1)}m\ot (b_0)S^{-1}(h^{(q+2)})\ot\dots\ot(a_q)S^{-1}(h^{(2)}).$$
 \end{lemma}
 \begin{proposition}
 For any $p\ge 0 $, the Hopf homology groups $\{H_p(\mathcal{H},M\ot B^{\ot(q+1)})\}_{q}$ is a cyclic module.
  \end{proposition}
  \begin{proof}
 This proposition  is true for any  cylindrical module with the same proof  given  in \cite{gj}.
  %Since $\{X_{p,q}\}_{p,q}$ is a cylindrical module and in particular the vertical
  % and horizontal operators commute among each other, it is follows that the operators are well-defined on 
   %the homology groups $\{ H_p(\mathcal{H},M\ot B^{\ot(q+1)})\}_{q\ge 0}$. Thus we only have to check the cyclic
    %condition ${t}^{q+1}=id$. 
 \end{proof}
 \begin{theorem}
There is a spectral sequence $E_{p,q}$ that converges to $HC^\mathcal{H}_\bullet(A,M)$.
 Its $E^2$-term is given by cyclic  homology groups,
 $$E^2_{p,q}=HC_q(H_p(\mathcal{H},M\ot B^{\ot(q+1)})).$$
\end{theorem} 
%\begin{corollary}
%Let $\mathcal{H}$ be  a Hopf algebra and $(\delta,\sigma)$ a modular pair in involution. Then 
% there is a spectral sequence that converges to $\widetilde(HC)^{(\delta,\sigma)}_{p+q}(\mathcal{H})$, 
% with $E^2$-terms  given by $E^2_{p,q}=H_p(\mathcal{H},HC_q(k))$.
%\end{corollary}
   It is easy to see that if $\mathcal{H}$ is semisimple,
    then the spectral sequence collapses at $E^1$-term. 
    Let $B^\natural_\mathcal{H}\ot M=H_0(\mathcal{H},M\ot b^{\ot(n+1)})$ 
  denote the first column of  $E^{1}$.
    
  \begin{proposition}
  Let $\mathcal{H}$ be a semisimple Hopf algebra, $B$ a right $\mathcal{H}$\nobreakdash-module algebra, and 
  $A=\mathcal{H}\# B$. Then we have $HC_n^\mathcal{H}(A,M)\cong HC_n(B^\natural_\mathcal{H}\ot M).$   
  \end{proposition}
  \begin{proof}
  Since $\mathcal{H}$ is semisimple, we have $E^1_{p,q}=H_p(\mathcal{H},M\ot B^{\ot(q+1)})=0$ for $p\ge 1$, 
  and the spectral
  sequence collapses. Then $E^2_{p,q}=0$ for  $p\ge 1$ and $E^2_{0,q}=HC_q(B^\natural_\mathcal{H}\ot M).$  
  \end{proof}

   %$$$$$$$$$$$$$$$$$$$$$$$$$$$$$$$$$$$$$$$$$$$$$$$$$$$$$$$$$$$$$$$$$$$$$$$$$$$$$$$$$$$$$$$$$$$          
 %%%%%%%%%%%%%%%%%%%%%%%%%%%%%%%%%%%%%%%%%%%%%%%%%%%%%%%%%%%%%%%%%%%%%%%%%%%%%%%%%%%%%%%%%%%%%%
 %%%%%%%%%%%%%%%%%%%%%%%%%%%%%%%%%%%%%%%%%%%%%%%%%%%%%%%%%%%%%%%%%%%%%%%%%%%%%%%%%%%%%%%%%%%%%%         
\section{ Invariant Cyclic Cohomology  of Hopf Cotriples }
In this section we define the  invariant cyclic cohomology of Hopf module coalgebras. One example is the 
Connes-Moscovici cyclic homology  of a Hopf algebra with a modular pair in involution 
in the sense of \cite{achm99} which turns  out to be the invariant cyclic cohomology  of the coalgebra $\mathcal{H}$. 
This  is implicit in \cite{achm98} and explicitly done in \cite{cr} for $\sigma=1$. We go, however, beyond this 
( fundamental) example and define  a cocyclic module for  any Hopf cotriple $(C,\mathcal{H},V)$ consisting  of  
an $\mathcal{H}$-module coalgebra $C$, an $\mathcal{H}$-comodule $V$ and a compatible character $\delta$ on
 $\mathcal{H}$. Some of the results that appear in Section \ref{sec3}, e.g. Morita invariance and the spectral 
  sequence in Subsection \ref{subsection3.4} can be formulated and proved in the context of Hopf cotriples.
   For the sake  of brevity we decide not to include these results here.     
\subsection{ Definition and basic properties}
\begin{definition}
By a left { Hopf cotriple} we mean a triple $(C,\mathcal{H},V)$ where $\mathcal{H}$ is a 
Hopf algebra, $C$ is a left $\mathcal{H}$\nobreakdash-module coalgebra and $V$ is 
a left  $\mathcal{H}$\nobreakdash-comodule.
\end{definition}
\begin{example}
Let $\mathcal{H}$ be a Hopf algebra and $V$ 
 a left  $\mathcal{H}$\nobreakdash-comodule. Then $(C,\mathcal{H},V)$ is a left Hopf cotriple,
 where $C=\mathcal{H}$ is the underlying coalgebra of $\mathcal{H}$, with  $\mathcal{H}$ acting on 
 $C$ via multiplication. In particular for $V=k_\sigma$ and $\mathcal{H}$ coacting on $k=k_\sigma$  via a grouplike 
 element $\sigma\in \mathcal{H}$, we have a Hopf cotriple $(\mathcal{H},\mathcal{H},k_\sigma)$. This is 
 the Hopf cotriple that is relevant to Connes-Moscovici theory \cite{achm98,achm99,achm00}. 
\end{example}
\begin{example}{(Trivial cotriples).} Let $C$ be a 
coalgebra, $\mathcal{H}=k$ and $V=k$. Then $(C,k,k)$ is a Hopf cotriple.
\end{example}

Given a Hopf cotriple $(C,\mathcal{H},V)$, let $C^n(C,V)=V\otimes C^{\ot(n+1)}$.
 We define cosimplicial and cyclic operators on $\{ C^n(C,V)\}_n=\{V\otimes C^{\ot(n+1)} \}_n$ by

\begin{eqnarray*}\label{f1}
\delta_i(v\otimes c_0 \otimes c_1 \otimes \dots \otimes c_n)&=& v\otimes c_0 \otimes\dots 
\otimes  c_i^{(1)}\otimes c_i^{(2)}\otimes c_n ~~~ 0 \leq i \leq n\\
\delta_{n+1}(v\otimes c_0 \otimes c_1 \otimes \dots \otimes c_n)&=& v^{({0})}\otimes c_0^{(2)}\otimes c_1 
\otimes \dots \otimes c_n  \otimes v^{({-1})}  c_0^{(1)}\\
\sigma_i(v\otimes c_0 \otimes c_1 \otimes \dots \otimes c_n)&=&v\otimes c_0 \otimes \dots 
c_i \otimes \varepsilon(c_{i+1})\otimes 
\dots\otimes c_n ~~~0\leq i \leq n-1\\
\tau(v\otimes c_0 \otimes c_1 \otimes \dots \otimes c_n)  &=&v^{({0})}\otimes c_1 \otimes c_2 \otimes 
\dots \otimes c_n \otimes v^{({-1})}c_0.
\end{eqnarray*}
\begin{proposition}
Endowed with the above  operators, $\{C^n(C,V)\}_{n}$ is a paracocyclic module.
\end{proposition}
  
  We have a diagonal $\mathcal{H}$-action on $C^n(C,V)$, defined by 
  $$h(v\ot c_o\ot\dots \ot c_n)=v\ot h^{(1)}c_0\ot h^{(2)}c_1\ot \dots \ot h^{(n+1)}c_n. $$
   It is easy to see that $C^n(C,V)$ is  an $\mathcal{H}$\nobreakdash-module.
   
   To define the space of coinvariants, we fix a character  of $\mathcal{H}$, say  $\delta$. Let
   $$C_\mathcal{H}^n(C,V)=\frac{ C^n(C,V)}{\text{span}
   \{hm-\delta(h)m\mid m\in C^n(C,V),\; h\in\mathcal{H}\}}$$
 be the space of coinvariants of $C^n(C,V)$ under the action of $\mathcal{H}$ and with respect
 to $\delta$.
 Our first task is to find sufficient  conditions under which  $\{C_\mathcal{H}^n(C,V)\}_n$ 
 is a cocyclic module. 
 
   Let us recall the twisted antipode $\widetilde{S}:\mathcal{H}\rightarrow \mathcal{H}$, where 
  $\widetilde{S}(h)=\delta(h^{(1)})S(h^{(2)})$, from ~\cite{achm98}.
  We define  the $V$-{\it twisted antipode}
  \begin{center}
  $\widetilde{S}_V:V\otimes\mathcal{H}\longrightarrow V\otimes\mathcal{H}$\\
  \shoveleft {by }\\
   $\widetilde{S}_V(v\otimes h)=v^{({0})}\otimes S^{-1}(v^{({-1})})\widetilde{S}(h). $
    \end{center}
    The following lemmas are  very useful in our theory.
  \begin{lemma}
  The twisted antipode and the $V$-twisted antipode are invertible.
    \end{lemma}
    \begin{proof}
    Define 
    \begin{align*}
    &&{\widetilde{S}}^{-1}:\mathcal{H}\longrightarrow \mathcal{H}, 
    &&&~~~~~{\widetilde{S}_V}^{-1}:V\ot \mathcal{H}\longrightarrow V\ot \mathcal{H}\\
    &\text{by} \\ 
   &&{\widetilde{S}}^{-1}(h)=\delta(h^{(2)})S^{-1}(h^{(1)}), 
    &&&~~~~~{\widetilde{S}_V}^{-1}(v\otimes h)=v^{(0)}\otimes {\widetilde{S}}^{-1}(h)S^{-1}(v^{(-1)}).
    \end{align*}
    One can check that $\widetilde{S}\circ{\widetilde{S}}^{-1}={\widetilde{S}}^{-1}
    \circ\widetilde{S}=id_\mathcal{H}$
    , and $\widetilde{S}_V\circ{\widetilde{S}_V}^{-1}={\widetilde{S}_V}^{-1}
    \circ\widetilde{S}_V=id_{V\otimes\mathcal{H}}.$
    \end{proof}
    \begin{lemma}{(\cite{cr})}
    Let $M$ and $N$ be two left $\mathcal{H}$\nobreakdash-modules. Then for all $h\in\mathcal{H}$, 
    $m\in M$, and $n\in N$ we have 
    $$hm\otimes n=m\otimes \widetilde{S}(h)n \quad \text{in}\quad (M\otimes N)_{\mathcal{H}}, $$ 
    where $-_{\mathcal{H}}$  is the functor of coinvariants with respect  to $\delta$. 
    \end{lemma} 
    \begin{definition}
    We call the pair $(\delta,V)$ a { comatched pair} if 
    $$v^{({0})}\delta(v^{({-1})})=v \;\;\;\text{for all}\; v\in V.$$ 
     We call the comatched pair $(\delta, V)$ a { comatched pair in involution} if
     $$(\widetilde{S}_V)^2=id_{V\otimes\mathcal{H}}.$$ 
    \end{definition}
    \begin{definition}
    Let $\delta$ be a character of $\mathcal{H}$. A Hopf cotriple 
    $(C,\mathcal{H},V)$ is called {$\delta$}-{ compatible} if $(\delta,V)$ is a comatched 
    pair in involution.
    \end{definition}
    \begin{lemma}
    If $(\delta,V)$ is a comatched pair,  then for any  $v\in V$ we have 
    $$v=v^{(0)}\delta(v^{(-1)})=v^{(0)}\delta(S(v^{(-1)}))=v^{(0)}\delta(S^{-1}(v^{(-1)})).$$
    \end{lemma}
    \begin{proof}
    The first equality is the definition of comatched pair. We prove the second one. Indeed\\ 
    $~~~~~~~~~~~~~~~~~~~~~~~~~~v^{(0)}\delta(S(v^{(-1)}))=v^{(0)}\delta(v^{(-1)})\delta
    (S(v^{(-2)}))=v^{(0)}\epsilon(v^{(-1)})=v.$\\
    The last equality can be proved in  a similar way. 
    \end{proof}
    Now we can prove the main result of this section. 
    \begin{theorem}{\label{theorem4.10}}
    Let $(C,\mathcal{H},V)$ be a $\delta$-compatible Hopf cotriple. Then $\{C^n_\mathcal{H}(C,V)\}_n$
     is a cocyclic module. 
     
    \end{theorem}
    \begin{proof}
    We should first  show that the cosimplicial and  cyclic operators on
     $\{C_\mathcal{H}^n(C,V)\}_n$ are well defined. In the following we check 
      this just for  the cyclic operator and the last coface. The rest is easy to check.
    
    Let $h\in\mathcal{H}$ and $v\otimes c_0\otimes \dots\otimes c_n\in C_\mathcal{H}^n(C,V)$. We  prove that\\ 
    $\tau(h(v\otimes c_0\otimes \dots\otimes c_n))=\delta(h)\tau(v\otimes c_0\otimes \dots\otimes c_n)$ 
    in the coinvariant space.  \\
    Indeed,
    \\
     $\tau(h(v\otimes c_0\otimes \dots\otimes c_n))=\tau(v\otimes h^{(1)}c_0\otimes \dots\otimes h^{(n+1)}c_n)$\\
     
     $=v^{({0})}\otimes h^{(2)}c_1\otimes \dots\otimes h^{(n+1)}c_n\otimes v^{({-1})}h^{(1)}c_0$\\
     
     $=h^{(2)}(v^{({0})}\otimes c_1\otimes \dots\otimes c_n)\otimes v^{({-1})}h^{(1)}c_0$\\
     
        $=\widetilde{S}^{-1}(S^{-1}(v^{({-2})}))\widetilde{S}^{-2}(h^{(2)})S^{-1}(v^{({-1})})
        (v^{({0})}\otimes c_1\otimes \dots\otimes c_n)\otimes v^{({-3})}h^{(1)}c_0$\\
        \\       
           $=v^{({0})}\otimes c_1\otimes\dots\otimes c_n\otimes\widetilde{S}
          (\widetilde{S}^{-1}(S^{-1}(v^{({-2})})))\widetilde{S}^{-2}(h^{(2)})S^{-1}(v^{({-1})})
         v^{({-3})}h^{(1)}c_0$\\
        \\
        $=v^{({0})}\otimes c_1\otimes \dots\otimes c_n\otimes
         \widetilde{S}(S^{-1}(v^{({-1})}))\widetilde{S}^{-1}(h^{(2)})
         S^{-1}(v^{({-2})})(v^{({-3})})h^{(1)}c_0$\\
        \\
        $=\delta(h)v^{({0})}\delta(S^{-1}v^{({-1})})\otimes c_1\otimes \dots\otimes
          c_n\otimes v^{({-2})} c_0$\\
         \\
         $=\delta(h)v^{({0})}\otimes c_1\otimes \dots\otimes c_n\otimes v^{({-1})} c_0$\\
         \\
         $=\tau(\delta(h)(v\otimes c_0\otimes \dots\otimes c_n)).$
         
         Now let us prove that the operator $\delta_{n+1}$ is also  well defined on  
         $C^n_\mathcal{H}(C,V)$. We have 
         \begin{align*}
         \delta_{n+1}(h(v\otimes c_0\otimes \dots\otimes c_n))&=
         \delta_{n+1}(v\otimes h^{(1)}c_0\otimes \dots\otimes h^{(n+1)}c_n)\\
         &=v^{({0})}\otimes h^{(2)}c_0^{(2)}\otimes \dots\otimes 
         h^{(n+1)}c_n\otimes v^{({-1})}h^{(1)}c_0^{(1)},
         \end{align*}
         and by the same steps as above this is equal 
         to $\delta_{n+1}(\delta(h)(v\otimes c_0\otimes \dots\otimes c_n)).$

         So $\{ C^n_\mathcal{H}(C,V)\}_n$ is at least a paracocyclic module.
         To prove that  it is cocyclic module
         we  check  that $\tau^{n+1}=id$.
         \begin{eqnarray*}
         \tau^{n+1}(v\otimes c_0\otimes \dots\otimes c_n)\hspace{-0.5cm}&&=v^{({0})}\otimes 
         (v^{({-n-1})})c_0\otimes \dots\otimes (v^{({-1})})c_n \\
         &&=v^{({-1})}(v^{({0})}\otimes c_0\otimes \dots\otimes c_n)\\
         &&=\delta(v^{({-1})})v^{({0})}\otimes c_0\otimes \dots\otimes c_n\\
         &&=v\otimes c_0\otimes \dots\otimes c_n.
         \end{eqnarray*}
    \end{proof}
    \begin{example}
    Let $\mathcal{H}=V=k$ Then  for any coalgebra $C$ one has $\{C_k^n(C,k)\}_n$ is 
    the natural cyclic module, $C_\natural$, of the coalgebra $C$.
    \end{example}
    \begin{example}\label{lemma4.12}
    Let $\mathcal{H}$ be a Hopf algebra,  $C=\mathcal{H}$, and $V=k_\sigma$. The Hopf cotriple 
    $(\mathcal{H},\mathcal{H},k_\sigma)$ is $\delta$-compatible if and only if $(\delta,\sigma)$
     is a modular pair in involution 
    in the sense of ~\cite{achm99}. In this case $\{ C^n_\mathcal{H}(C,V)\}_n$ is isomorphic 
    to the Connes-Moscovici cocyclic module $\mathcal{H}^{(\delta,\sigma)}_\natural$.
    \end{example}
    \begin{lemma}{\cite{cr}}
    Let $V$ be a left $\mathcal{H}$-module and $\delta$ a character of $\mathcal{H}$  then 
    $$ (\mathcal{H}\ot V)_\mathcal{H}\cong V,$$ 
    where $\mathcal{H}$ acts diagonally  on $\mathcal{H}\ot V$. 
    \end{lemma}
    \begin{example}
    We  generalize  the previous example. Let $\mathcal{H}$ be a Hopf algebra, $V$ be a $\mathcal{H}-$comodule,  
    and assume that $C=\mathcal{H}$. Then if $(\mathcal{H},\mathcal{H},V)$ is $\delta$-compatible we 
    can use  Lemma \ref{lemma4.12}  to simplify the cocyclic module
     $\{C^n_\mathcal{H}(\mathcal{H},V)\}_{n} $  
    to $\{V\otimes\mathcal{H}^{\otimes n}\}_{n}$ with the following  operators 
    \begin{align*}
    &\delta_0(v\otimes h_1\otimes \dots \otimes h_n)=v\otimes\ 1\otimes h_1\otimes \dots \otimes h_n\\
    &\delta_i(v\otimes h_1\otimes \dots \otimes h_n)=v\otimes h_1\otimes \dots 
    \Delta(h_i)\otimes \dots \otimes h_n, &&\hspace{-55pt} 1\le i\le n\\
    &\delta_{n+1}(v\otimes h_1\otimes \dots \otimes h_n)=v^{({0})}\otimes
     h_1\otimes \dots \otimes h_n\otimes v^{({-1})}&&\\
    &\sigma_i(v\otimes h_1\otimes \dots \otimes h_n)=v\otimes h_1\otimes \dots 
     \otimes h_{i-1}\otimes \epsilon( h_i)\otimes h_n, &&\hspace{-55pt} 0\le i\le n-1\\
    &\tau(v\otimes h_1\otimes \dots \otimes h_n)=v^{({0})}\otimes{S}(h_1^{(n)})
     h_2\otimes\dots\otimes S(h_1^{(2)})h_n\otimes \widetilde{S}(h_1^{(1)})v^{({-1})}.
      \end{align*} 
      \end{example}
   By a {\it cotrace } in $\mathcal{H}$ we mean an 
element $t\in \mathcal{H}$ with $t^{(1)}\otimes t^{(2)}=t^{(2)}\otimes t^{(1)}$.
 Let $\delta$ be a character of $\mathcal{H}$.
 \begin{definition}
 An element $t\in\mathcal{H}$ is called a left $\delta$-{\it integral} if for all $g\in \mathcal{H}$,
 $$tg=\delta(g)t.$$
 \end{definition}
 Obviously, for $\delta=\epsilon$, an $\epsilon$-integral is 
 simply a (left) integral in $\mathcal{H}$. 
 \begin{example}
 Let $G$ be a finite group and $\delta $ a character of $G$. Then one can check that
 $$t=\sum_{g\in G}\frac{g}{\delta(g)}$$
 is a cotracial $\delta$-integral in $\mathcal{H}=kG$.
 \end{example}
 \begin{proposition}
 If $\mathcal{H}$ admits  a cotracial $\delta$-integral $t$ with $\delta(t)$ invertible in $k$, then
 $HC_\mathcal{H}^\bullet(C,k)$ is a direct summand in $HC^\bullet(C)$.
 \end{proposition}
 \begin{proof}
 Given such an element $t$, we define a right inverse for the projection map 
 $$\pi: C^\bullet(C)\longrightarrow C_\mathcal{H}^\bullet(C,k)$$
 by 
 \begin{center}
 $\gamma:C_\mathcal{H}^\bullet(C,k) \longrightarrow    C^\bullet(C)  $\\
 $\gamma(c_0\otimes c_1\otimes \dots \otimes c_n)=
 t^{(1)}c_0\otimes t^{(2)}c_1\otimes \dots \otimes t^{(n+1)}c_0.$
 \end{center}
 One can check that $\gamma $ is a  well-defined cyclic module map and $\pi\gamma =\delta(t)id$.
 \end{proof}
 \subsection{The commutative case}
 In this subsection we show that if $\mathcal{H}$ is a commutative Hopf 
 algebra, then the invariant  cyclic cohomology groups of $(\mathcal{H},\mathcal{H},V)$
 decompose as 
\begin{equation}\label{e5}
 HC^n_\mathcal{H}(\mathcal{H},V)=\bigoplus_{i\ge 0}H^{n-2i}(\mathcal{H},V),
\end{equation}
 where on the right hand side we have the Hopf cohomology groups of $\mathcal{H}$ with
  coefficients in $V$ as defined below.
   
    Consider the functor of coinvariants $$V\mapsto V^{\text{co}\mathcal{H}}=
  \{v\in V\mid \rho(v)=1\otimes v\}$$  
  from the category of left $\mathcal{H}$\nobreakdash-comodules to the category of $k$-modules.
   Since $V^{\text{co}\mathcal{H}}=Hom_{\text{comod}}(k,V),$ 
  this functor is left exact. The corresponding right  derived functors 
  are denoted by $H^i(\mathcal{H},V)$, and they appear on (\ref{e5}). Alternatively,
   these groups can be calculated from the complex 
   $$V\overset{d_0}{\rightarrow}V\otimes\mathcal{H}\overset{d_1}
   {\rightarrow}V\otimes\mathcal{H}^{\ot{2}}\overset{d_2}{\rightarrow}\dots$$ 
   with the differentials $d_n$,  $n\ge 1$,  given by
   \begin{multline*}
   d_n(v\otimes h_1\otimes h_2\otimes \dots\otimes h_n)=v\otimes 1\otimes 
   h_1\otimes h_2\otimes \dots\otimes h_n+\\
   +\sum_{i=1}^n(-1)^iv\otimes h_1\otimes h_2\otimes
    \dots\otimes h_i^{(1)}\otimes h_i^{(2)}\otimes \dots \otimes h_n+
	(-1)^{n+1}v^{(0)}\otimes h_1\otimes h_2\otimes \dots\otimes h_n\otimes v^{(-1)}
  \end{multline*}
 and $d_0(v)=v\otimes 1-v^{(0)}\otimes v^{(-1)}.$

For any cosimplicial module $M=\{M^n\}_n$,  its ``path space" $EM$ is defined by $(EM)^n=M^{n+1}$ 
with all operators shifted by 
$1$. It is a cosimplicial module. In particular, for $\{C_\mathcal{H}^n(\mathcal{H},V)\}_n$ 
we obtain  the cosimplicial module 
$\{EC_\mathcal{H}^n(\mathcal{H},V)\}_n.$
\begin{proposition}If $\mathcal{H}$ is commutative, then 
$\{EC_\mathcal{H}^n(\mathcal{H},V)\}_n$ has a cocyclic module 
structure. 
\end{proposition}  
\begin{proof}
Define a cyclic action
\begin{center}
$t:EC_\mathcal{H}^n(\mathcal{H},V)\longrightarrow EC_\mathcal{H}^n(\mathcal{H},V)$\\
\makebox[8in][l]{by}\\
$t(v\otimes h_0\otimes h_1\otimes \dots \otimes h_n)=v^{(0)}\otimes h_0^{(1)}S(h_1^{(n)})h_2\otimes
\dots \otimes h_0^{(n-1)}S(h_1^{(2)})h_n\otimes h_0^{(n)}S(h_1^{(1)})v^{(-1)}.$
\end{center}
One can check that all the axioms for cocyclic modules are satisfied.
\end{proof}
Now $\mathcal{H}$ has a  coaction on $EC_\mathcal{H}^n(\mathcal{H},V)$ by 
\begin{center}
$\rho:EC_\mathcal{H}^n(\mathcal{H},V)\longrightarrow\mathcal{H}\ot EC_\mathcal{H}^n(\mathcal{H},V)$
$\rho(v\ot h_0\ot h_1\ot\dots\ot h_n)=h_0^{(1)}\ot v\ot h_0^{(2)}\ot h_1\ot\dots\ot h_n.$
    \end{center}
   One can see that $(EC_\mathcal{H}^n(\mathcal{H},V))^{\text{co}\mathcal{H}}\cong C_\mathcal{H}^n(\mathcal{H},V)$.
    By the same method that we used in the proof of Proposition \ref{hasan2}, we can prove the following 
     proposition. 
    \begin{proposition}\label{hasan3}
    Let $\mathcal{H}$ be a commutative Hopf algebra. Then the invariant  cyclic cohomology  of the 
    $\epsilon$-compatible Hopf cotriple $(\mathcal{H},\mathcal{H},V)$ are given by
    $$HC_\mathcal{H}^n(\mathcal{H},V)=\bigoplus_{i\ge 0}H^{n-2i}(\mathcal{H},V).$$
    \end{proposition}
    
    Let $G$ be a complex affine algebraic group and $G\times M\longrightarrow M$ a linear action of $G$
    on a finite dimensional complex vector space $M$. Then $V=M^\ast$ is a comodule 
    over $\mathcal{H}=\mathbb{C}\lbrack G \rbrack$. The cohomology groups $H^i(\mathcal{H},V)$
     are easily seen to be isomorphic  to (algebraic) group cohomology $H^i_{alg}(G,M)$, 
     where  cochains $f:G\times G\times \dots \times G\rightarrow M$ are  assumed to be  algebraic functions.
      It therefore follows from Proposition \ref{hasan3}, that  
      $$HC_\mathcal{H}^n(\mathcal{H},M)=\bigoplus_{i\ge 0}H^{n-2i}(G,M).$$

    \end{document}